\documentclass[graybox]{svmult}

\usepackage{type1cm}        % activate if the above 3 fonts are
\usepackage{makeidx}         % allows index generation
\usepackage{graphicx}        % standard LaTeX graphics tool
\usepackage{multicol}        % used for the two-column index
\usepackage[bottom]{footmisc}% places footnotes at page bottom

\usepackage{newtxtext}       % 
\usepackage{newtxmath}       % selects Times Roman as basic font

\input epsf

\def\arcsinh{{\rm arcsinh}}
\def\calC{{{\cal C}}}
\def\calL{{{\cal L}}}
\def\bigO{{\cal O}}
\def\RR{\mathbb R}

\def\ph{{\rm ph}}

\def\sk{\sum_{k=0}^\infty\,}

\def\Frac#1#2{\frac
{
 {\raise.6ex
 \hbox{$\displaystyle#1$}}
}
{
 {\lower.6ex
 \hbox{$\displaystyle#2$}}
 }
}
\def\Ai{{{\rm Ai}}}
\newcommand{\dsp}{\displaystyle}
\newcommand{\arccosh}{{\rm arccosh}}

\makeindex             % used for the subject index
\begin{document}

\title*{Asymptotic computation of classical orthogonal polynomials}
\author{A. Gil, J. Segura and N. M. Temme}
\institute{A. Gil \at Departamento de Matem\'atica Aplicada y CC. de la Computaci\'on.
ETSI Caminos. Universidad de Cantabria. 39005-Santander, Spain. \email{amparo.gil@unican.es}
\and J. Segura \at Departamento de Matem\'aticas, Estad\'{\i}stica y 
        Computaci\'on, Universidad de Cantabria, 39005 Santander, Spain. \email{segurajj@unican.es}
\and N.M. Temme \at IAA, 1825 BD 25, Alkmaar, The Netherlands 
(Former address: CWI, Science Park 123, 1098 XG Amsterdam,  The Netherlands). \email{nicot@cwi.nl}}
\maketitle

\abstract{
The classical orthogonal polynomials (Hermite, Laguerre and Jacobi) are involved in a vast 
number of applications in physics and engineering.
When large degrees $n$ are needed, the use of recursion to compute the polynomials
is not a good strategy for computation and a more efficient approach, such as the use of asymptotic expansions,
is recommended. In this paper, we give an overview of the asymptotic expansions considered in \cite{Gil:2017:ELP}
for computing Laguerre polynomials $L^{(\alpha)}_n(x)$ for bounded values of the parameter $\alpha$.
Additionally,
we show examples 
of the computational performance of an asymptotic expansion for $L^{(\alpha)}_n(x)$ 
valid for large values of $\alpha$ and $n$. This expansion was 
used in \cite{Gil:2018:GHL} as starting point for obtaining asymptotic approximations to the zeros.
Finally, we analyze the expansions 
considered in \cite{Gil:2018:ELB}, \cite{Gil:2019:ELP} and \cite{Gil:2019:NIJ}
to compute the Jacobi polynomials for large degrees $n$.
}

\section{Introduction}

As it is well known, Laguerre $L^{(\alpha)}_n(x)$ and Jacobi  $P_{n}^{(\alpha,\beta)}(x)$ polynomials 
present monotonic and oscillatory regimes depending on the parameter 
values. Sharp bounds limiting the oscillatory region for these polynomials are given in 
\cite{Dimitrov:2010:BEZ}. In the particular case of Jacobi polynomials, these bounds are

\begin{equation}\label{intro2}
\Frac{B-4(n-1)\sqrt{C}}{A} \le  x_{n,k}(\alpha,\beta) \le \Frac{B+4(n-1)\sqrt{C}}{A},
\end{equation}
where $x_{n,k}(\alpha,\beta)$ are the zeros of the Jacobi polynomials and

\begin{equation}\label{intro3}
\begin{array}{lcl}
A&=&(2n+\alpha+\beta)(n(2n+\alpha+\beta)+2(\alpha+\beta+2)),\\
B&=&(\beta-\alpha)((\alpha+\beta+6)n+2(\alpha+\beta)),\\
C&=&n^2(n+\alpha+\beta+1)^2+(\alpha+1)(\beta+1)(n^2+(\alpha+\beta+4)n+2(\alpha+\beta)).
\end{array}
\end{equation}

\begin{figure}
\epsfxsize=14cm \epsfbox{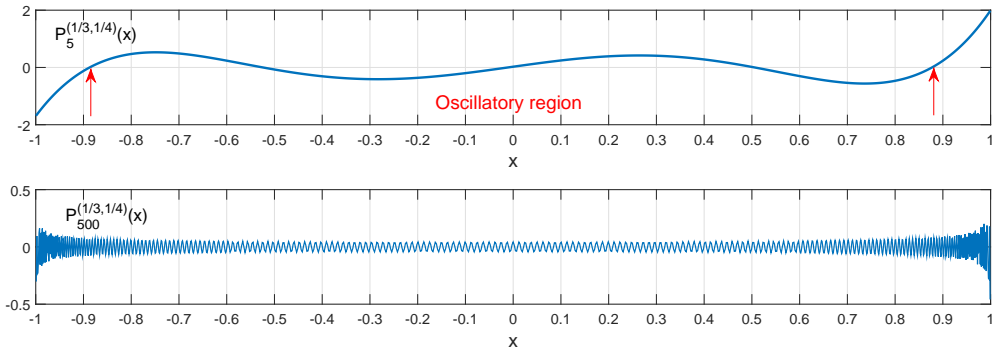}
\caption{
\label{fig:fig01} Two plots of Jacobi polynomials ($P^{(1/3,1/4)}_5(x)$ and  $P^{(1/3,1/4)}_{500}(x)$)
showing the oscillatory behaviour of the polynomials.
}
\end{figure}

An example of the oscillatory behaviour of the Jacobi polynomials for two different values of $n$
is given in Figure~\ref{fig:fig01}.

On the other hand, the classical orthogonal polynomials satisfy (also well-known) three-term recurrence relations

\begin{equation}\label{eq:Intro04}
p_{n+1}(x)=(A_n x+B_n)p_n(x)-C_np_{n-1},
\end{equation}
where, for Laguerre polynomials $L^{(\alpha)}_n(x)$, the coefficients are given by 

\begin{equation}\label{eq:Intro05}
A_{n}=\dfrac{-1}{n+1},\, B_{n}=\dfrac{2n+\alpha+1}{n+1},\,
C_{n}=\dfrac{n+\alpha}{n+1},
\end{equation}
and for Jacobi polynomials $P_{n}^{(\alpha,\beta)}(x)$,

\begin{equation}\label{eq:Intro05}
\begin{array}{@{}r@{\;}c@{\;}l@{}}
A_{n}&=&\dfrac{(2n+\alpha+\beta+1)(2n+\alpha+\beta+2)}{2(n+1)(n+\alpha+\beta+1)},\\[8pt]
B_{n}&=&\dfrac{(\alpha^{2}-\beta^{2})(2n+\alpha+\beta+1)}{2(n+1)(n+\alpha+\beta+%
1)(2n+\alpha+\beta)},\\[8pt]
C_{n}&=&\dfrac{(n+\alpha)(n+\beta)(2n+\alpha+\beta+2)}{(n+1)(n+\alpha+\beta+1)(2%
n+\alpha+\beta)}.
\end{array}
\end{equation}

By using the Perron-Kreuser theorem \cite{Kreuser:1914:UDV} it is easy to show that the three-term recurrence relations for the
classical orthogonal polynomials are not badly conditioned in the sense that there is no solution which dominates
exponentially over the rest. 
Therefore the recurrence relations can be used in the forward 
direction, for example, with starting values  $L^{(\alpha)}_0(x)=1$
and $L^{(\alpha)}_1(x)=1+\alpha-x$ for Laguerre polynomials, or 
$P_{0}^{(\alpha,\beta)}(x)=1$ and $P_{1}^{(\alpha,\beta)}(x)=\tfrac12(\alpha-\beta)+\tfrac12(\alpha+\beta+2)x$,
to compute the Jacobi polynomials
when $n$ is small/moderate.
 However, when large degrees $n$ are needed, the use of recursion to compute the polynomials
is not a good strategy for computation and a more efficient approach, such as the use of asymptotic expansions,
is recommended.

The use of asymptotic expansions
 in the computation of Hermite polynomials has already been discussed in \cite{Gil:2006:PCF}, 
since Hermite polynomials are a particular case of the parabolic cylinder function $U(a,x)$.
In this paper, we consider the asymptotic computation of Laguerre and Jacobi polynomials. 
Other recent references dealing with the numerical evaluation of these polynomials 
using asymptotic methods are, for example,
\cite{Dea:2016:Jac} and \cite{Huy:2018:Lag}.

 We first give an overview of the asymptotic expansions considered in \cite{Gil:2017:ELP}
for computing Laguerre polynomials $L^{(\alpha)}_n(x)$ for bounded values of the parameter $\alpha$.
Additionally,
we show examples 
of the computational performance of an asymptotic expansion for $L^{(\alpha)}_n(x)$ given in \cite{Gil:2018:GHL}
valid for large values of $\alpha$ and $n$. 
We also analyze the expansions 
considered in \cite{Gil:2018:ELB}, \cite{Gil:2019:ELP} and \cite{Gil:2019:NIJ}
to compute the Jacobi polynomials for large degrees $n$.
These expansions have been used as starting point to obtain asymptotic approximations 
to the nodes  and/or weights of the Gauss-Jacobi (G-J) quadrature.  
In the asymptotic expansions for large degree $n$ described in
 \cite{Gil:2019:NIJ}, we assume that $\alpha$ and $\beta$ should be of order ${\cal O}(1)$ for large values of $n$. 
It is interesting to note that these approximations 
can be used as standalone
methods for the non-iterative computation ($15-16$ digits accuracy) of the nodes  and weights of G-J quadratures of high degree
($n \ge 100$). Alternatively, they can be used as accurate starting points of methods for solving nonlinear equations. This is
the approach followed, for example, in \cite{Hale:2013:SJSC} in which simple asymptotic approximations for the nodes are iteratively
refined by Newton's method.

\section{Laguerre polynomials}

The computation of Laguerre polynomials $L^{(\alpha)}_n(z)$ for bounded values of the parameter $\alpha$ was considered in 
\cite{Gil:2017:ELP}. In that paper, we analyzed the computational performance of three types of asymptotic expansions valid for large
values of $n$: an expansion in terms of Airy functions and two expansions in terms of Bessel functions (one for
small values of the variable $z$ and other ones in which larger values are allowed.) It is possible to obtain asymptotic expansions in terms of elementary functions starting from a Bessel expansions, as done in \cite{Dea:2013:SRA}; however, the Bessel expansions appear to have a larger range of validity. 

We next summarize results given in \cite{Gil:2017:ELP} and \cite{Gil:2018:GHL}; more details about the expansions can 
be found in \cite{Temme:2015:AMI}.

\begin{description}
\item{\bf a)  An expansion in terms of Airy functions } 

The representation used is

\begin{equation}\label{eq:lagairy01}
L_n^{(\alpha)}(\nu x)=(-1)^n
\frac{ e^{\frac12\nu x}\chi(\zeta)}{2^\alpha\nu^{\frac13}}\left(\Ai\left(\nu^{2/3} \zeta\right)A(\zeta)
+\nu^{-\frac43}\Ai^{\prime}\left(\nu^{2/3}\zeta\right)
B(\zeta)\right)
\end{equation}
with  expansions
\begin{equation}\label{eq:lagairy02}
A(\zeta)\sim\sum_{j=0}^\infty\frac{\alpha_{2j}}{\nu^{2j}},\quad B(\zeta)\sim\sum_{j=0}^\infty\frac{\beta_{2j+1}}{\nu^{2j}},\quad \nu \to\infty.
\end{equation}

Here
\begin{equation}\label{eq:lagairy03}
\nu=4\kappa,\quad \kappa=n+\tfrac12(\alpha+1), \quad 
\chi(\zeta)=2^{\frac12}x^{-\frac14-\frac12\alpha}\left(\frac{\zeta}{x-1}\right)^{\frac14},
\end{equation}
and

\begin{equation}
 \label{eq:lagairy04} 
\begin{array}{ll}
\dsp{\tfrac23(-\zeta)^{\frac32}=\tfrac12\left(\arccos\sqrt{{x}}-\sqrt{{x-x^2}}\right)}  & \mbox{if $ 0<x \le 1$,}\\
\dsp{\tfrac23\zeta^{\frac32}=\tfrac12\left(\sqrt{{x^2-x}}-\arccosh\sqrt{{x}}\right)} & \mbox{if $ x\ge 1$.}
\end{array}
\end{equation}

The expansions in \eqref{eq:lagairy02} hold uniformly for bounded $\alpha$ and  $x\in(x_0,\infty)$, where $x_0\in(0,1)$, a fixed number. We observe that $x=1$ is a turning point, where the argument of the Airy functions vanishes (see \eqref{eq:lagairy04}). Left of $x=1$, where $\zeta <0$, the zeros of $L_n^{(\alpha)}(\nu x)$ occur.

The first few coefficients in the expansions (\ref{eq:lagairy02}) are $\alpha_0=1$ and

\begin{description}
\item{For $ 0<x \le 1$}
$$
\begin{array}{lcl}
\alpha_2&=&-\Frac{1}{1152}\left(4608u^8v^6\mu^3+2016v^6u^8\mu^2+2304u^8v^6\mu^4 \right.\\
            &&    -672v^6u^4\mu^2-192v^6\mu u^6+2880v^6u^6\mu^2+2688v^3u^7\mu^2\\
             &&   -288v^6\mu u^8+3072u^6v^6\mu^3-672v^6\mu u^4+2688v^3\mu u^7+\\
              &&   924v^6u^2+558v^6u^4-180u^6v^6-135v^6u^8+504u^7v^3+336v^3u^5\\
              &&   \left.+280u^3v^3-7280u^6+385v^6\right)/(v^6u^6),\\
\beta_1&=&\Frac{1}{12}\left(48v^3u^4\mu^2+9v^3u^4+48v^3\mu u^4-20u^3+6v^3u^2+5v^3\right)/(u^3v^4),
\end{array}
$$
where $u= \sqrt{1/x-1}$,  $v= 2\sqrt{-\zeta}$ and $\mu= (\alpha-1)/2$.

\item {For $x>1$}
$$
\begin{array}{lcl}
\alpha_2&=& \Frac{1}{1152}\left(4608u^8v^6\mu^3+2016v^6u^8\mu^2+\right.\\
            && 2304u^8v^6\mu^4-672v^6u^4\mu^2+192v^6\mu u^6 \\
           &&-2880v^6u^6\mu^2-336v^3u^5-288v^6\mu u^8-3072u^6v^6\mu^3\\
           &&-672v^6\mu u^4+280u^3v^3-924v^6u^2+558v^6u^4+180u^6v^6\\
           &&-135v^6u^8+2688v^3u^7\mu^2+504u^7v^3+2688v^3\mu u^7\\
           &&\left.-7280u^6+385v^6\right)/(u^6v^6),\\
\beta_1&=& \Frac{1}{12}\left(48v^3u^4\mu^2+9v^3u^4+48v^3\mu u^4-20u^3-6v^3u^2+5v^3\right)/(u^3v^4),
\end{array}
$$
where $u=\sqrt{1-1/x}$, $v= 2\sqrt{\zeta}$ and $\mu= (\alpha-1)/2$.
\end{description}

A detailed explanation of the method used to obtain the coefficients of Airy expansions can be found
in \cite[\S23.3]{Temme:2015:AMI}. The expansion in terms of Airy function was derived in  \cite{Frenzen:1988:UAE}.

\item{\bf b) A simple Bessel-type expansion}

The starting point to obtain this expansion is the well-known
relation between Laguerre polynomials and the Kummer function ${}_1F_1(a;c;x)$

\begin{equation}\label{eq:lagsim01}
L_n^{(\alpha)}(z)=
\left(
\begin{array}{c}
n+\alpha\\
n
\end{array}
\right)
{}_1 F_1
\left(
\begin{array}{c}
-n\\
\alpha+1
\end{array}
;z
\right)
\,.
\end{equation}

Then, we use an expansion for  ${}_1F_1(a;c;x)$ for large negative values of $a$; 
see  \cite[\S10.3.4]{Temme:2015:AMI}

\begin{equation}\label{eq:lagsim02}
\begin{array}{ll}
\dsp{\frac{1}{\Gamma(c)}
{}_1 F_1
\left(
\begin{array}{c}
-a\\
c
\end{array}
;z
\right)
\sim\left(\frac{z}{a}\right)^{\frac12(1-c)}\frac{\Gamma(1+a)e^{\frac12z}}{\Gamma(a+c)}\ \times}\\[8pt]
\quad\quad\dsp{\left(J_{c-1}\left(2\sqrt{az}\right)\sk \frac{a_k(z)}{(-a)^k}-\sqrt{\frac{z}{a}}J_{c}\left(2\sqrt{az}\right)\sk \frac{b_k(z)}{(-a)^k}\right).}
\end{array}
\end{equation}
This expansion of ${}_1F_1(-a;c;z)$ is valid for bounded values of $z$ and $c$, with $a\to\infty$ inside the sector  $-\pi+\delta\le\ph\,a\le\pi-\delta$.
Using this, we obtain

\begin{equation}\label{eq:lagsim03}
\begin{array}{ll}
\dsp{L_n^{(\alpha)}(x)\sim\left(\frac{x}{n}\right)^{-\frac12\alpha} e^{\frac12x}\ \times}\\[8pt]
\quad\dsp{\left(J_{\alpha}\left(2\sqrt{nx}\right)\sk (-1)^k\frac{a_k(x)}{n^k}-\sqrt{\frac{x}{n}}J_{\alpha+1}\left(2\sqrt{nx}\right)\sk (-1)^k \frac{b_k(x)}{n^k}\right),}
\end{array}
\end{equation}
as $n\to\infty$. 

The coefficients $a_k(x)$ and $b_k(x)$ follow from the expansion of the function

\begin{equation}\label{eq:lagsim03}
f(x,s)=e^{xg(s)}\left(\frac{s}{1-e^{-s}}\right)^{\alpha+1},\quad g(s)=\frac{1}{s}-\frac{1}{e^s-1}-\frac12.
\end{equation}

The function $f$ is analytic in the strip $\vert\Im s\vert<2\pi$ and it can be expanded for  $\vert s\vert<2\pi$ into
\begin{equation}\label{eq:lagsim04}
f(x,s)=\sum_{k=0}^\infty c_k(x) s^k.
\end{equation}

The coefficients $a_k(x)$ and $b_k(x)$ are given in terms of the $c_k(x)$ coefficients: 

\begin{equation}\label{eq:lagsim07}
\begin{array}{@{}r@{\;}c@{\;}l@{}}
a_k(x) & = & \dsp{\sum_{m=0}^k 
\left(
\begin{array}{c}
k\\
m
\end{array}
\right)
(m+1-c)_{k-m}x^m c_{k+m}(x),}\\[8pt]
b_k(x) & = & \dsp{\sum_{m=0}^k 
\left(
\begin{array}{c}
k\\
m
\end{array}
\right)
(m+2-c)_{k-m}x^m c_{k+m+1}(x),}
\end{array}
\end{equation}
$k=0,1,2,\ldots$.

The coefficients $c_k(x)$  are combinations of Bernoulli numbers and Bernoulli polynomials, the first ones being (with $c=\alpha+1$)
\begin{equation}\label{eq:lagsim06}
\begin{array}{@{}r@{\;}c@{\;}l@{}}
c_0(x)&=&1,\quad c_1(x)=\frac{1}{12}\left(6c-x\right),\quad \\[8pt]
c_2(x)&=&\frac{1}{288}\left(-12c+36c^2-12xc+x^2\right),\\[8pt]
c_3(x)&=&\frac{1}{51840}\left(-5x^3 + 90x^2c +(-540c^2 + 
 180c+72)x +1080c^2(c-1)\right),
\end{array}
\end{equation}
and the first relations are
\begin{equation}\label{eq:lagsim08}
\begin{array}{ll}
a_0(x)= c_0(x)=1,\quad b_0(x)= c_1(x),\\[8pt]
a_1(x)= (1-c)c_1(x)+xc_2(x),\quad b_1(x)= (2-c)c_2(x)+xc_3(x),\\[8pt]
a_2(x)= (c^2-3c+2)c_2(x)+(4x-2xc)c_3(x)+x^2c_4(x),\\[8pt]
b_2(x)=  (c^2-5c+6)c_3(x)+(6x-2xc)c_4(x)+x^2c_5(x),
\end{array}
\end{equation}
again with $c=\alpha+1$.

\item{\bf c) A not so simple expansion in terms of Bessel functions} 

To understand how this representation and the coefficients of the expansion are obtained, we start with the standard form
\begin{equation}\label{eq:Bessalgcoef01}
F_\zeta(\nu)=\frac1{2\pi i}\int_{{\calC}}e^{\nu\left(u-\zeta/u\right)}f(u)\,\frac{du}{u^{\alpha+1}},
\end{equation}
where the contour ${\calC}$ starts at $-\infty$ with  $\ph\,u=-\pi$, encircles the origin anti-clockwise, 
and returns to $-\infty$ with $\ph\, u=\pi$. The function $f(u)$ is assumed to be analytic in a neighborhood of $\calC$,
 and  in particular in a domain that contains the saddle points $\pm ib$, where $b=\sqrt{\zeta}$. 

When $f$ is replaced by unity, we obtain the Bessel function:
\begin{equation}\label{eq:Bessalgcoef02}
F_\zeta(\nu)=\zeta^{-\frac12\alpha} J_\alpha\left(2\nu\sqrt{\zeta}\right). 
\end{equation}

For Laguerre polynomials, we have the integral representation
\begin{equation}\label{eq:lagbess001}
L_n^{(\alpha)}(z)= \frac{1}{2\pi i}\int_\calL (1-t)^{-\alpha-1}e^{-tz/(1-t)}\,\frac{dt}{t^{n+1}},
\end{equation}
where $\calL$ is a circle around the origin with radius less than unity.

After substituting $t=e^{-s}$, we obtain 
\begin{equation}\label{eq:lagbess002}
e^{-\nu x}L_n^{(\alpha)}(2\nu x)=\frac{2^{-\alpha}}{2\pi i}\int_{-\infty}^{(0+)}e^{\nu h(s,x)}\left(\frac{\sinh s}{s}\right)^{-\alpha-1}\,\frac{ds}{s^{\alpha+1}},
\end{equation}
where $\nu=2n+\alpha+1$ and $h(s,x)=s-x\coth s$. The contour is the same as in \eqref{eq:Bessalgcoef01}.  
Using the transformation $h(s,x)=u-\zeta/u$, we obtain the integral representation in the standard form 
of \eqref{eq:Bessalgcoef01}

\begin{equation}\label{eq:lagbessa}
2^{\alpha}e^{-\nu x}L_n^{(\alpha)}(2\nu x)=\frac{1}{2\pi i}\int_{-\infty}^{(0+)}e^{\nu(u-\zeta/u)}f(u)\,\frac{du}{u^{\alpha+1}},
\end{equation}
where
\begin{equation}\label{eq:lagbessb}
f(u)=\left(\frac{u}{\sinh s}\right)^{\alpha+1}\frac{ds}{du}.
\end{equation}

We now consider the following recursive scheme
\begin{equation}\label{eq:Bessalgcoef03}
\begin{array}{@{}r@{\;}c@{\;}l@{}}
f_j(u)&=&A_j(\zeta)+B_j(\zeta)/u+\left(1+b^2/u^2\right)g_j(u),\\[8pt]
 \dsp{f_{j+1}(u)}&=&\dsp{g_j^\prime(u)-\frac{\alpha+1}{u}g_j(u),}\\[8pt]
A_j(\zeta)&=&\dsp{\frac{f_j(ib)+f_j(-ib)}{2}},\quad \dsp{B_j(\zeta)=i\frac{f_j(ib)-f_j(-ib)}{2b},}
\end{array}
\end{equation}
with $f_0(u)=f(u)$, the coefficient function. 

Using  this scheme and integration by parts in  \eqref{eq:lagbessa}, we obtain the following representation

\begin{equation}\label{eq:lagbess01}
L_n^{(\alpha)}(2\nu x)= \frac{e^{\nu x}\chi(\zeta)}{2^\alpha\zeta^{\frac12\alpha}}\left(J_\alpha\bigl(2\nu \sqrt{\zeta}\bigr)
A(\zeta)-
\frac{1}{\sqrt{\zeta}}J_{\alpha+1}\bigl(2\nu  \sqrt{\zeta}\bigr)B(\zeta)\right),
\end{equation}
with expansions
\begin{equation}\label{eq:lagbess02}
A(\zeta)\sim\sum_{j=0}^\infty\frac{A_{2j}(\zeta)}{\nu^{2j}},\quad 
B(\zeta)\sim\sum_{j=0}^\infty\frac{B_{2j+1}(\zeta)}{\nu^{2j+1}},\quad \nu\to\infty.
\end{equation}
Here,
\begin{equation}\label{eq:lagbess03}
\nu=2n+\alpha+1,\quad \chi(\zeta)=(1-x)^{-\frac14}\left(\frac{\zeta}{x}\right)^{\frac12\alpha+\frac14},\quad x<1,
\end{equation}
with $\zeta$ given by

\begin{equation}\label{eq:lagbess04}
\begin{array}{ll}
\dsp{\sqrt{-\zeta}=\tfrac12\left(\sqrt{{x^2-x}}+\arcsinh\sqrt{{-x}}\right)}, & \quad\mbox{if \quad$x\le0$,}\\[8pt]
\dsp{\sqrt{\zeta}=\tfrac12\left(\sqrt{{x-x^2}}+\arcsin\sqrt{{x}}\right)}, & \quad \mbox{if \quad$ 0\le x<1$.}
\end{array}
\end{equation}

This expansion in terms of Bessel functions was derived in  \cite{Frenzen:1988:UAE}.  It is uniformly valid for $x\le 1-\delta$, where $\delta\in(0,1)$  is a fixed number. Recall that the expansion in terms of Airy functions is valid around the turning point $x=1$ and up to infinity, but not for small positive $x$. As shown in  
\cite{Frenzen:1988:UAE}, the  Bessel-type expansion (left of $x=1$) and the Airy-type expansion are valid in overlapping $x$-domains of $\RR$. For $x<0$ the $J$-Bessel functions can be written in terms of modified $I$-Bessel functions.

The coefficients $A_j(\zeta)$ and $B_j(\zeta)$ in \eqref{eq:lagbess02} can all be expressed in terms of the derivatives $f^{(k)}(\pm ib)$
 of $f(u)$ (see  \eqref{eq:lagbessb}) at the saddle points $\pm ib$.  We have to proceed as follows:

We expand the functions $f_j(u)$ \eqref{eq:Bessalgcoef03} in two-point Taylor expansions
\begin{equation}\label{eq:Bessalgcoef05}
f_j(u)= \sum_{k=0}^\infty  C_k^{(j)} (u^2-b^2)^k+u\sum_{k=0}^\infty  D_k^{(j)} (u^2-b^2)^k.
\end{equation}

Using the recursive scheme for the functions $f_j(u)$   given in \eqref{eq:Bessalgcoef03}, 
we derive the following scheme for the coefficients
\begin{equation}\label{eq:Bessalgcoef06}
\begin{array}{ll}
\dsp{C_k^{(j+1)}=(2k-\alpha)D_{k}^{(j)}+b^2(\alpha-4k-2)D_{k+1}^{(j)}+2(k+1)b^4D_{k+2}^{(j)},}\\[8pt]
\dsp{D_k^{(j+1)}=(2k+1-\alpha)C_{k+1}^{(j)}-2(k+1)b^2C_{k+2}^{(j)},}
\end{array}
\end{equation}
for $j,k=0,1,2,\ldots$, and the coefficients $A_j$ and $B_j$ follow from 
\begin{equation}\label{eq:Bessalgcoef07}
A_j(\zeta)=C_0^{(j)},\quad B_j(\zeta)=-b^2D_0^{(j)},\quad j \ge 0.
\end{equation}

In the present case of the Laguerre polynomials the functions $f_{2j}$ are even and $f_{2j+1}$ are odd, and we have $A_{2j+1}(\zeta)=0$ and $B_{2j}(\zeta)=0$. A few non--vanishing coefficients are
\begin{equation}\label{eq:Bessalgcoef08}
\begin{array}{@{}r@{\;}c@{\;}l@{}}
A_0(\zeta)&=&f(ib),\\[8pt]
B_1(\zeta)&=&-\frac{1}{4}b\bigl((2\alpha-1)if^{(1)}(ib)+bf^{(2)}(ib)\bigr),\\[8pt]
A_2(\zeta)&=&-\frac{1}{32b}\bigl(3i(4^2\alpha-1)f^{(1)}(ib)-
(3-16\alpha+4\alpha^2)bf^{(2)}(ib)\ \\[8pt]
&&+2i(2\alpha-3)b^2f^{(3)}(ib)b^2+b^3f^{(4)}(ib)\bigr),\\[8pt]
B_3(\zeta)&=&-\frac{1}{384b}\bigl(3(4\alpha^2-1)(2\alpha-3)(if^{(1)}(ib)+bf^{(2)}(ib))\ \\[8pt]
&&+2i(\alpha-7)(2\alpha-1)(2\alpha-3)b^2f^{(3)}(ib)\ \\[8pt]
&&+3(19-20\alpha+4\alpha^2))b^3f^{(4)}(ib)\\[8pt]
&&-3i(2\alpha-5)b^4f^{(5)}(ib)-b^5f^{(6)}(ib)\bigr)\\.
\end{array}
\end{equation}

To have $A_0(\zeta)=1$  we have to scale all $A$ and $B$-coefficients with respect to $A_0(\zeta)=\chi(\zeta)$.
As an example, using this  normalization and the explicit expression for the first two derivatives $f^{(k)}(ib)$
 of $f(u)$, we obtain for $B_1(\zeta)$ the following expression

\begin{equation}
B_1(\zeta)=\dsp{ \frac{1}{48\xi}\left(5\xi^4b+6\xi^2b+3\xi+12\alpha^2(b-\xi)-3b\right),}
\end{equation}
where $\dsp{\xi=\sqrt{\frac{x}{1-x}}}$.

In the resulting algorithm for Laguerre polynomials described in \cite{Gil:2017:ELP},
the two Bessel-type expansions were used in the oscillatory region of the functions.
The use of the simple Bessel-type expansion was restricted, as expected, to values of 
the argument of the Laguerre
polynomials very close to the origin. The second Bessel-type expansion was used in the part
of the oscillatory region where the Airy expansion failed (computational schemes for Airy functions
are more efficient than the corresponding schemes for  Bessel functions.) Therefore, the Airy expansion 
was used in the
monotonic region and in 
part  of the oscillatory region of the Laguerre polynomials $L^{(\alpha)}_n(x)$, with $\alpha$ small.

Recall that we have observed (see below \eqref{eq:lagbess04}) that the Airy-type and Bessel-type expansions in 
\eqref{eq:lagairy01} and \eqref{eq:lagbess01}, respectively, are valid in overlapping $x$-intervals of $\RR$.

\item{\bf d) An expansion in terms of Bessel functions for large $\alpha$ and $n$} 

When $\alpha$ is large, an expansion suitable for the oscillatory region of the Laguerre polynomials and
valid for large values of $n$, is the expansion given in \cite{Gil:2018:GHL}

\begin{equation}\label{eq:LagBesla01}
\begin{array}{@{}r@{\;}c@{\;}l@{}}
L_n^{(\alpha)}(4\kappa x)&=&\dsp{e^{-\kappa A}\chi(b)\left(\frac{b}{2\kappa x}\right)^{\alpha}\frac{\Gamma(n+\alpha+1)}{n!}
\ \times }\\[8pt]
&&
\quad\quad
\dsp{\bigl(J_\alpha(4\kappa b)A(b)-2bJ_\alpha^\prime(4\kappa x)B(b)\bigr),}
\end{array}
\end{equation}
with expansions
\begin{equation}\label{eq:LagBesla02}
A(b)\sim\sum_{k=0}^\infty \frac{A_k(b)}{\kappa^k}, \quad
B(b)\sim\sum_{k=0}^\infty \frac{B_k(b)}{\kappa^k},\quad \kappa\to\infty,
\end{equation}
where
\begin{equation}\label{eq:LagBesla03}
\kappa=n+\tfrac12(\alpha+1),\quad \chi(b)=\left(\frac{4b^2-\tau^2}{4x-4x^2-\tau^2}\right)^{\frac14},\quad \tau=\frac{\alpha}{2\kappa}. \end{equation}
We have $\tau <1$. The quantity  $b$ is a function of $x$ and follows from the relation
\begin{equation}\label{eq:LagBesla04}
\begin{array}{@{}r@{\;}c@{\;}l@{}}
&&\dsp{2W-2\tau\arctan\frac{W}{\tau}}=\\[8pt]
&&\quad\quad\quad
\dsp{2R-\arcsin\frac{1-2x}{\sqrt{1-\tau^2}}-\tau\arcsin\frac{x-\frac12\tau^2}{x\sqrt{1-\tau^2}}+\tfrac12\pi(1-\tau),}
\end{array}
\end{equation}
where 
\begin{equation}\label{eq:LagBesla05}
R=\tfrac12\sqrt{4x-4x^2-\tau^2}=\sqrt{(x_2-x)(x-x_1)},\quad W=\sqrt{4b^2-\tau^2},
\end{equation}
and
\begin{equation}\label{eq:LagBesla06}
x_1=\tfrac12\left(1-\sqrt{1-\tau^2}\right), \quad 
x_2=\tfrac12\left(1+\sqrt{1-\tau^2}\right).
\end{equation}

The relation in \eqref{eq:LagBesla04} can be used for $x\in[x_1,x_2]$, in which case $b\ge\frac12\tau$. In this interval the
Laguerre polynomial $L_n^{(\alpha)}(4\kappa x)$ oscillates.  The endpoints of this interval can be seen as two turning points, where in fact Airy functions can be used. We assume that $\tau\le\tau_0$, where $\tau_0$ is a fixed number in the interval $(0,1)$. When $\tau\to1$, the points $x_1$ and $x_2$, which can be considered as turning points, become the same, and we expect expansions in terms of parabolic cylinder function (in the present case these are Hermite polynomials).

The expansion is valid for $x\in [0,(1-\delta)x_2]$, with $\delta$ a small positive number.  For $x < x_1$, the relation in \eqref{eq:LagBesla04} should be modified, because the arguments of the square roots of $W$ and $R$ become negative. 

At the turning point $x=x_1$ we have $b=\frac12\tau$, and the argument of the Bessel functions in \eqref{eq:LagBesla01} becomes $4\kappa b=2\kappa\tau=\alpha$, where we have used the definition of $\tau$ in \eqref{eq:LagBesla03}. We see that argument and order of the Bessel functions become equal, and because we assume that $\alpha$ is large, this explains that at the turning point $x_1$ we can expect Airy-type behaviour. The Airy-type approximation around $x_1$ cannot be valid near $x=0$, but the Bessel function can handle this.

The first coefficients of the expansions in \eqref{eq:LagBesla02} are 
\begin{equation}\label{eq:LagBesla07}
\begin{array}{@{}r@{\;}c@{\;}l@{}}
A_0(b)&=&1,\quad B_0(b)=0,\\[8pt]
A_1&=&\dsp{\frac{\tau}{24(\tau^2-1)}}, \quad B_1=\dsp{\frac{P R^3+QW^3}{192R^3W^4(\tau^2-1)}},\\[8pt]
&&\quad P= 4(2\tau^2+12b^2)(1-\tau^2),\\[8pt]
&&\quad Q= 2\tau^4-12x^2\tau^2-\tau^2-8x^3+24x^2-6x,\\[8pt]
B_2(b&)=&A_1(b)B_1(b).
\end{array}
\end{equation}

An example of the computational performance of the expansion  (\ref{eq:LagBesla01}) is given
in Figure~\ref{fig:fig02}. For comparison, we show the
results for two different values of $\alpha$ ($\alpha=198.5$, $\alpha=20.6$) and
a fixed value of  $n$ ($n=340$). The coefficients $A_k, \,B_k,\,k=0,1,2$ in (\ref{eq:LagBesla02}) have been considered
in the computations.
The expansion has been tested 
by considering the relation given in Eq.~(18.9.13) of  \cite{Koornwinder:2010:OPS},
written in the form

\begin{equation}\label{eq:comp1}
\epsilon=\left|\Frac{L_{n-1}^{(\alpha+1)}(z)+L_n^{(\alpha)}(z)}{L_n^{(\alpha+1)}(z)}-1\right|,
\end{equation}
where $z=4\kappa x$ is the argument of the Laguerre polynomial in the expansion (\ref{eq:LagBesla01}). 
As can be seen, an accuracy better than $10^{-10}$ can be obtained in the first half of the $x$-interval 
for the two values of $\alpha$ tested.

\end{description}
\begin{figure}
%\begin{center}
\epsfxsize=14cm \epsfbox{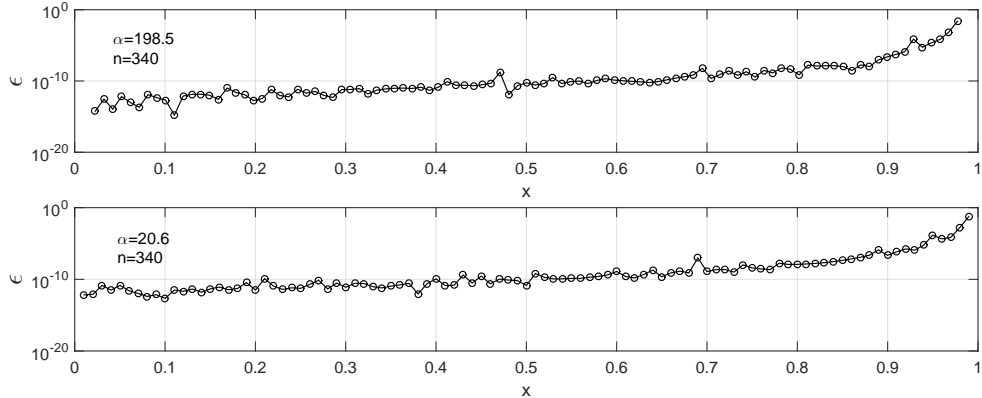}
\caption{\label{fig:fig02}  Test of the performance of the expansion (\ref{eq:LagBesla01}) for large values of $n$ and $\alpha$. 
The $\epsilon$-value in Eq. (\ref{eq:comp1}) is plotted. }
%\end{center}
\end{figure}

\section{Jacobi polynomials}\label{sec:Jacexpan}

Next we review of some useful asymptotic expansions for computing Jacobi polynomials.

\subsection{An expansion in terms of elementary functions}\label{sec:jacelem}

In \cite{Gil:2019:NIJ} we considered the following representation

\begin{equation}\label{eq:Jacelem01}
P_n^{(\alpha,\beta)}(\cos\theta)=\frac{G_\kappa(\alpha,\beta)}{\sqrt{\pi \kappa}}\frac{\cos \chi \,U(x)-\sin\chi \,V(x)}
{\sin^{\alpha+\frac12}\frac12\theta\,\cos^{\beta+\frac12}\frac12\theta},
\end{equation}
where
\begin{equation}\label{eq:Jacelem02}
x=\cos\theta,\quad \chi=\kappa\theta-\left(\tfrac12\alpha+\tfrac14\right)\pi,\quad \kappa = n+\tfrac12(\alpha+\beta+1).
\end{equation}

The expansion is valid for $x\in[-1+\delta,1-\delta]$, where $\delta$ is a small positive number. 

The functions $U(x)$ and $V(x)$ have the expansions

\begin{equation}\label{eq:Jacelem03}
U(x)\sim\sum_{k=0}^\infty \frac{u_{2k}(x)}{\kappa^{2k}},\quad V(x)\sim\sum_{k=0}^\infty \frac{v_{2k+1}(x)}{\kappa^{2k+1}}.
\end{equation}

The first coefficients are
\begin{equation}\label{eq:Jacelem04}
\begin{array}{@{}r@{\;}c@{\;}l@{}}
u_0(x)&=&1, \quad  v_1(x)= \dsp{\frac{2\alpha^2-2\beta^2+(2\alpha^2+2\beta^2-1)x}{8\sin\theta} } ,\\[8pt]
u_2(x)&=&\dsp{\frac{1}{384\sin^2\theta}}\bigl(
12(5-2\alpha^2-2\beta^2)(\alpha^2-\beta^2)x\ +\\[8pt]
&&\quad
4(-3(\alpha^2-\beta^2)^2+3(\alpha^2+\beta^2)-6+4\alpha(\alpha^2-1+3\beta^2)\ +\\[8pt]
&&\quad
(-12(\alpha^2+\beta^2)(\alpha^2+\beta^2-1)-16\alpha(\alpha^2-1+3\beta^2)-3)x^2\bigr).
\end{array}
\end{equation}

The function $G_\kappa(\alpha,\beta)$ is given by
 
\begin{equation}\label{eq:Jacelem05}
G_\kappa(\alpha,\beta)=\frac{\Gamma(n+\alpha+1)}{n!\,\kappa^\alpha}=\frac{\Gamma\left(\kappa+\frac12(\alpha-\beta+1)\right)}{\Gamma\left(\kappa-\frac12(\alpha+\beta-1)\right)\,\kappa^\alpha},
\end{equation}
with expansion in negative powers of $\left(\kappa-\frac12\beta\right)^2$ 

 \begin{equation}\label{eq:Jacelem6}
G_\kappa(\alpha,\beta)\sim (w/\kappa)^{\alpha}\sum_{m=0}^\infty \frac{C_m(\rho) (-\alpha)_{2m}}{w^{2m}},
\end{equation}
where
 \begin{equation}\label{eq:Jacelem7}
w=\kappa-\tfrac12\beta, \quad \rho=\tfrac12(\alpha+1),
\end{equation}
and the first few $C_m(\rho)$ coefficients are
 \begin{equation}\label{eq:Jacelem8}
\begin{array}{l}
C_0(\rho)=1, \quad C_1(\rho)=-\tfrac1{12}\rho, \quad C_2(\rho)=\tfrac1{1440}\left(5\rho+1\right),\\
C_3(\rho)= -\rho \Frac{4+21\rho+35\rho^2}{362880},\\[8pt]
C_4(\rho)= \rho\Frac{18+101\rho+210\rho^2+175\rho^3}{87091200}.
\end{array}
\end{equation}

\subsection{An expansion in terms of Bessel functions}\label{sec:jacBess}

Another expansion, this time in terms of Bessel functions, considered in \cite{Gil:2019:NIJ}
is

\begin{equation}\label{eq:JacBes01}
\begin{array}{@{}r@{\;}c@{\;}l@{}}
P_n^{(\alpha,\beta)}(\cos\theta)&=&\dsp{\frac{G_\kappa(\alpha,\beta)}{\sin^{\alpha}\frac12\theta\,
\cos^{\beta}\frac12\theta}\sqrt{\frac{\theta}{\sin\theta}}\,W(\theta)}, \\[8pt]
W(\theta)&=&\dsp{J_{\alpha}(\kappa\theta)\,S(\theta)+\frac{1}{\kappa}J_{\alpha+1}(\kappa\theta)\,T(\theta),}
\end{array}
\end{equation}
where  $G_\kappa(\alpha,\beta)$ is defined in \eqref{eq:Jacelem05},  and
\begin{equation}\label{eq:JacBes02}
\kappa=n+\tfrac12(\alpha+\beta+1).
\end{equation}

This expansion holds uniformly with respect to $\theta \in[0,\pi-\delta]$, where $\delta$ is a small positive number.
When an asymptotic expansion is needed for $\theta$ near $\pi$ (or for $x$ near $-1$), we can use the symmetry relation
\begin{equation}\label{eq:JacBessym}
P_n^{(\alpha,\beta)}(-x)=(-1)^nP_n^{(\beta,\alpha)}(x).
\end{equation}

The functions $S(\theta)$ and $T(\theta)$ have the expansions
\begin{equation}\label{eq:JacBes03}
S(\theta)\sim \sum_{k=0}^\infty\frac{S_k(\theta)}{\kappa^{2k}},\quad
T(\theta)\sim \sum_{k=0}^\infty\frac{T_k(\theta)}{\kappa^{2k}},\quad \kappa\to\infty,
\end{equation}
with $S_0(\theta)=A_0(\theta)=1$,  $T_0(\theta)=A_1(\theta)$, and for $k=1,2,3,\ldots$
\begin{equation}\label{eq:JacBes04}
\begin{array}{@{}r@{\;}c@{\;}l@{}}
S_k(\theta)&=&\dsp{-\frac{1}{\theta^{k-1}}\sum_{j=0}^{k-1}
\left( 
\begin{array}{c}
k-1\\
j
\end{array}
\right)
A_{j+k+1}(\theta)(-\theta)^j2^{k-1-j}(\alpha+2+j)_{k-j-1}},\\[8pt]
T_k(\theta)&=&\dsp{\frac{1}{\theta^{k}}\sum_{j=0}^{k} 
\left( 
\begin{array}{c}
k\\
j
\end{array}
\right)
A_{j+k+1}(\theta)(-\theta)^j2^{k-j}(\alpha+1+j)_{k-j}}.  
\end{array}
\end{equation}

The coefficients $A_k(\theta)$ are analytic functions for $0\le\theta<\pi$.  
The first ones are $A_0(\theta)=1$ and
\begin{equation}\label{eq:JacBes05}
A_1(\theta)=\frac{\left(4\alpha^2-1\right)(\sin\theta-\theta\cos\theta)+2\theta\left(\alpha^2-\beta^2\right)(\cos\theta-1)}{8\theta\sin\theta}.
\end{equation}

For small values of $\theta$, it is convenient to consider expansions of the coefficients  $A_k(\theta)$:

\begin{equation}\label{eq:JacBes10}
A_k(\theta)=\chi^k \theta^{k}\sum_{j=0}^\infty A _{jk}\theta^{2j},\quad \chi=\frac{\theta}{\sin\theta},
\end{equation}
where the series represent entire functions of $\theta$. The first few $A_{jk}$ are
\begin{equation}\label{eq:JacBes11}
\begin{array}{@{}r@{\;}c@{\;}l@{}}
A_{0,1}&=& \dsp{\tfrac{1}{24}\left(\alpha^2+3\beta^2-1\right), }\\[8pt]
A_{1,1}&=& \dsp{\tfrac{1}{480}\left(-3\alpha^2-5\beta^2+2\right), }\\[8pt]
A_{0,2}&=& \dsp{\tfrac{1}{5760}\left(-16\alpha-14\alpha^2-90\beta^2+5\alpha^4+4\alpha^3+45\beta^4+30\beta^2\alpha^2+60\beta^2\alpha+21\right). }
\end{array}
 \end{equation}

More details about the coefficients  $A_k(\theta)$  are given in \cite{Gil:2019:NIJ}.

As in the case of the Bessel-type expansions for Laguerre polynomials, for computing the representation (\ref{eq:JacBes01}),
an algorithm for evaluating the Bessel functions $J_{\nu}(z)$ is needed.
 We use a computational scheme based on the methods of approximation described in the Appendix.

\subsection{An expansion for large $\beta$ }\label{sec:jaclargeb}

In \cite{Gil:2018:ELB} we discussed two asymptotic approximations of Jacobi polynomials for large values of the $\beta$-parameter.
The expansions are in terms of Laguerre polynomials.
One of the two expansions given in \cite{Gil:2018:ELB} is the following:

\begin{equation}\label{eq:Jacbeta1}
\begin{array}{@{}r@{\;}c@{\;}l@{}}
\dsp{P_n^{(\alpha,\beta)}\left(1-\frac{2z}{b}\right)}
&=&(1-z/b)^n W(n,\alpha,z),\\[8pt]
W(n,\alpha,z)&=&L_{n}^{(\alpha)}(z)R(n,\alpha,\beta,z)+L_{n-1}^{(\alpha)}(z)S(n,\alpha,\beta,z),
\end{array}
\end{equation}
where $0<z<\beta$, $b=\beta+n$ and $L_{n}^{(\alpha)}(z)$, $L_{n-1}^{(\alpha)}(z)$ are Laguerre polynomials. $R(n,\alpha,\beta,z)$,
$S(n,\alpha,\beta,z)$ have the expansions 
\begin{equation}\label{eq:Jacbeta2}
R(n,\alpha,\beta,z)\sim\sum_{k=0}^\infty\frac{r_k(n,\alpha,z)}{b^k},\quad S(n,\alpha,\beta,z)\sim\sum_{k=0}^\infty\frac{s_k(n,\alpha,z)}{b^k}.
\end{equation}

The first few coefficients are

\begin{equation}\label{eq:Jacbeta3}
\begin{array}{@{}r@{\;}c@{\;}l@{}}
r_0(n,\alpha,z)&=&1,\quad  s_0(n,\alpha,z)=0,\\[8pt]
r_1(n,\alpha,z)&=&\frac12n(2z+\alpha+1),\\[8pt]
s_1(n,\alpha,z)&=&-\frac12(n+\alpha)(\alpha+z+1),\\[8pt]
r_2(n,\alpha,z)&=&-\frac{1}{24}\left(3\alpha^3+6\alpha^2z-6\alpha nz+3\alpha z^2\right.\\
     &&\left.-9nz^2+10\alpha^2+8\alpha z-4nz-12z^2+9\alpha+2\right),\\[8pt]
s_2(n,\alpha,z)&=& \frac{1}{24}\left(n+\alpha)(3\alpha^3+3\alpha^2z-6\alpha n z-3\alpha z^2-6nz^2\right.\\
    &&\left.-3z^3+10\alpha^2+\alpha z-4nz-11z^2+9\alpha-2z+2\right).
\end{array}
\end{equation}

In this expansion for the Jacobi polynomials $P_n^{(\alpha,\beta)}(x)$, the range of validity in the $x$-interval is very limited near $x=1$ 
(to reach the whole $x$-interval $(-1,1)$, $z$ should be of order $\bigO(\beta)$). 
For large values of the parameter $\alpha$, we can obtain a similar expansion valid near $x=-1$
 by using the symmetry relation \eqref{eq:JacBessym} of the Jacobi polynomials.

\subsection{An expansion for large $\alpha$ and $\beta$ }\label{sec:jaclargeab}

In  \cite{Gil:2019:ELP} an expansion for Jacobi polynomials in terms of elementary functions valid for large values
of the parameters $\alpha$ and $\beta$ has been given:

\begin{equation}\label{eq:Jacalphabeta1}
P_n^{(\alpha,\beta)}(x)=\frac{2^{\frac12(\alpha+\beta+1)}e^{-\kappa\psi}}
{\sqrt{\pi \kappa w(x)U(x)}}\left(\cos\left(\kappa\chi+\tfrac14\pi\right)P+\sin\left(\kappa\chi+\tfrac14\pi\right)Q\right),
\end{equation}
where

\begin{equation}\label{eq:Jacalphabeta2}
\begin{array}{@{}r@{\;}c@{\;}l@{}}
w(x)&=&(1-x)^\alpha(1+x)^\beta,\\[8pt]
\psi&=&-\frac12(1-\tau)\ln(1-\tau)-\frac12(1+\tau)\ln(1+\tau)\ +\\[8pt]
&&\frac12(1+\sigma)\ln(1+\sigma)+\frac12(1-\sigma)\ln(1-\sigma),\\[8pt]
\chi&=&\dsp{(\tau+1)\arctan\frac{U(x)}{1-x+\sigma+\tau}\ +}\\[8pt]
&&\dsp{(\tau-1)\arctan\frac{U(x)}{1+x+\sigma-\tau}
%-(1-\sigma)\arctan\frac{U(x)}{\tau+x\sigma}.}
+(1-\sigma)\,{\rm{atan}}2(-U(x),\tau+x\sigma).}\\[8pt]
U(x)&=&\sqrt{1-2\sigma\tau x-\tau^2-\sigma^2-x^2},
\end{array}
\end{equation}
and
\begin{equation}\label{eq:Jacalphabeta3}
\sigma=\frac{\alpha+\beta}{2\kappa},\quad \tau=\frac{\alpha-\beta}{2\kappa},\quad \kappa=n+\tfrac12(\alpha+\beta+1).
\end{equation}

The function ${\rm{atan}}2(y, x)$ in the third term of $\chi(x)$ in \eqref{eq:Jacalphabeta2} denotes the phase
$\in (-\pi,\,\pi]$ of the complex number $x + iy$. 

$P$ and $Q$ have expansions

\begin{equation}\label{eq:Jacalphabeta4}
P\sim \sum_{j=0}^\infty \frac{p_{j}}{\kappa^j},\quad Q\sim \sum_{j=0}^\infty \frac{q_{j}}{\kappa^j},
\end{equation}
and the first coefficients are $p_0=1$, $q_0=0$. 
 These expansions are valid for $x\in[x_-(1+\delta),x_+(1-\delta)]$, where $x_\pm$ are the turning points 

$$
x_\pm=-\sigma\tau\pm\sqrt{(1-\sigma^2)(1-\tau^2)},
$$
and $\delta$ is a fixed positive small number. 
For more details about the coefficients $p_j$ and $q_j$, see \cite{Gil:2019:ELP}.

\subsection{Computational performance}

The performance of the asymptotic expansions for the Jacobi polynomials has been tested 
by considering the relation given in Eq.(18.9.3) of  \cite{Koornwinder:2010:OPS}
written in one of the two following forms

\begin{equation}\label{eq:test1}
\Frac{P_{n}^{(\alpha,\,\beta-1)}(x)-P_n^{(\alpha-1,\,\beta)}(x)}{P_{n-1}^{(\alpha,\,\beta)}(x)}=1,
\end{equation}

\begin{equation}\label{eq:test2}
\Frac{P_{n}^{(\alpha-1,\,\beta)}(x)+P_{n-1}^{(\alpha,\,\beta)}(x)}{P_{n}^{(\alpha,\,\beta-1)}(x)}=1. 
\end{equation}

Note that the test \eqref{eq:test1} fails close to the zeros of  $P_{n-1}^{(\alpha,\,\beta)}(x)$; in this case, we can consider
the alternative test \eqref{eq:test2}: the zeros of $P_{n-1}^{(\alpha,\,\beta)}(x)$ and $P_{n}^{(\alpha,\,\beta-1)}(x)$ interlace according to 
Theorem 2 of \cite{seg:2008:int}, and therefore they 
can not vanish simultaneously. 

\begin{figure}
%\begin{center}
\epsfxsize=14cm \epsfbox{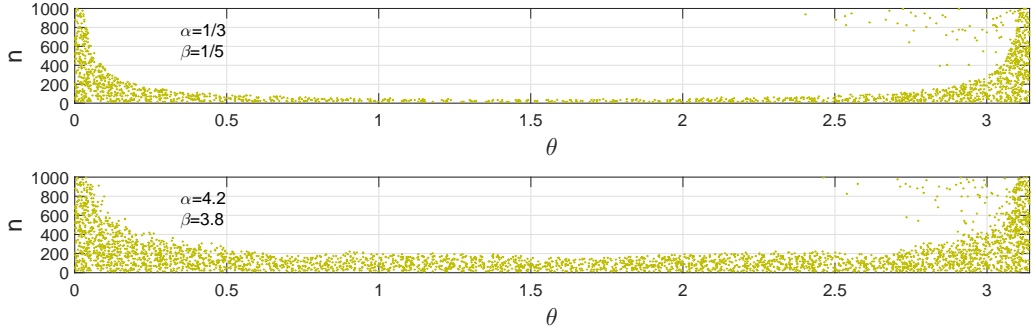}
\caption{
\label{fig:fig03} Performance of the asymptotic expansion (\ref{eq:Jacelem01}). 
The points where the error in the computation 
 is greater than
$10^{-12}$. 
Two different sets of values of the parameters $\alpha$ and $\beta$ 
have been considered in the tests.}
%\end{center}
\end{figure}

First, we test the performance of the asymptotic expansion (\ref{eq:Jacelem01}).
Figure~\ref{fig:fig03} shows the points where the error in the computation of the test is greater than
$10^{-12}$  for two different sets of values of the parameters $(\alpha, \,\beta)$.
The first four nonzero coefficients of the expansions (\ref{eq:Jacelem03})
have been considered in the computations. Random values 
have been generated in the parameter region $(n,\,\theta) \in (10,\,1000) \times (0,\,\pi)$.  
As expected, the expansion works well for values of $x=\cos (\theta)$ not close to $\pm 1$. Also,
 the parameters $\alpha,\,\beta$ should be of order $\bigO(1)$, and smaller values of these parameters give better results, as can be seen when comparing the
results for $P_{n}^{(1/3,\,1/5)}(\cos \theta)$ and  $P_{n}^{(4.2,\,3.8)}(\cos \theta)$.

\begin{figure}
%\begin{center}
\epsfxsize=14cm \epsfbox{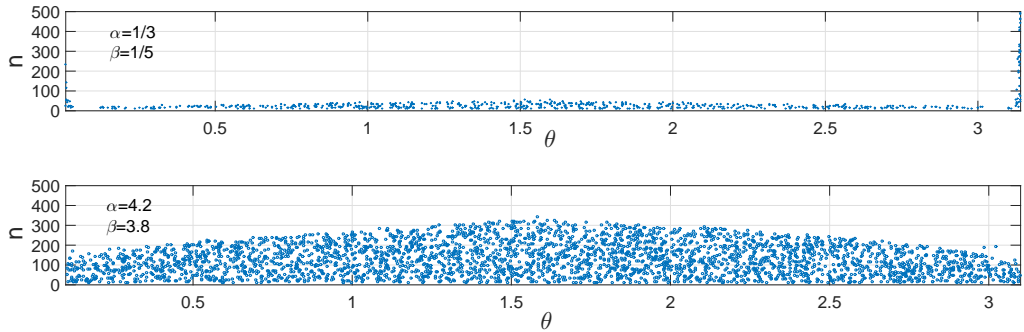}
\caption{
\label{fig:fig04}
 Tests of the performance of the asymptotic expansion
 (\ref{eq:JacBes01}) for two different sets of values of the parameters $\alpha$ and $\beta$ .
The points where the error in the computation  is greater than
$10^{-12}$  are shown. 
}
%\end{center}
\end{figure}

The performance of the asymptotic expansion in terms of Bessel functions (\ref{eq:JacBes01}) is illustrated in
Figure~\ref{fig:fig04}. As before, the points where the
error in the computation is greater than
$10^{-12}$  are shown. 
The asymptotic expansion
(\ref{eq:JacBes01}) has been computed by using the first three coefficients  $S_k(\theta),\,T_k(\theta),\,k=0,1,2$ in (\ref{eq:JacBes03}).
The results are as expected: better accuracy is obtained near $x=\pm 1$, although the expansion performs well 
even for small values of $x$ ($\theta$-values close to $\pi/2$) when $\alpha$ is small.  Also, for a fixed number of terms used in the expansion,
the accuracy obtained with the expansion decreases as the
parameter values $(\alpha,\,\beta)$ increase.

Regarding CPU times, the Bessel expansions, as expected, are slower than the
elementary expansions. Therefore, when both expansions are valid in the same parameter region, it is always convenient 
to use the expansion in terms of elementary functions.

\begin{figure}
%\begin{center}
\epsfxsize=14cm \epsfbox{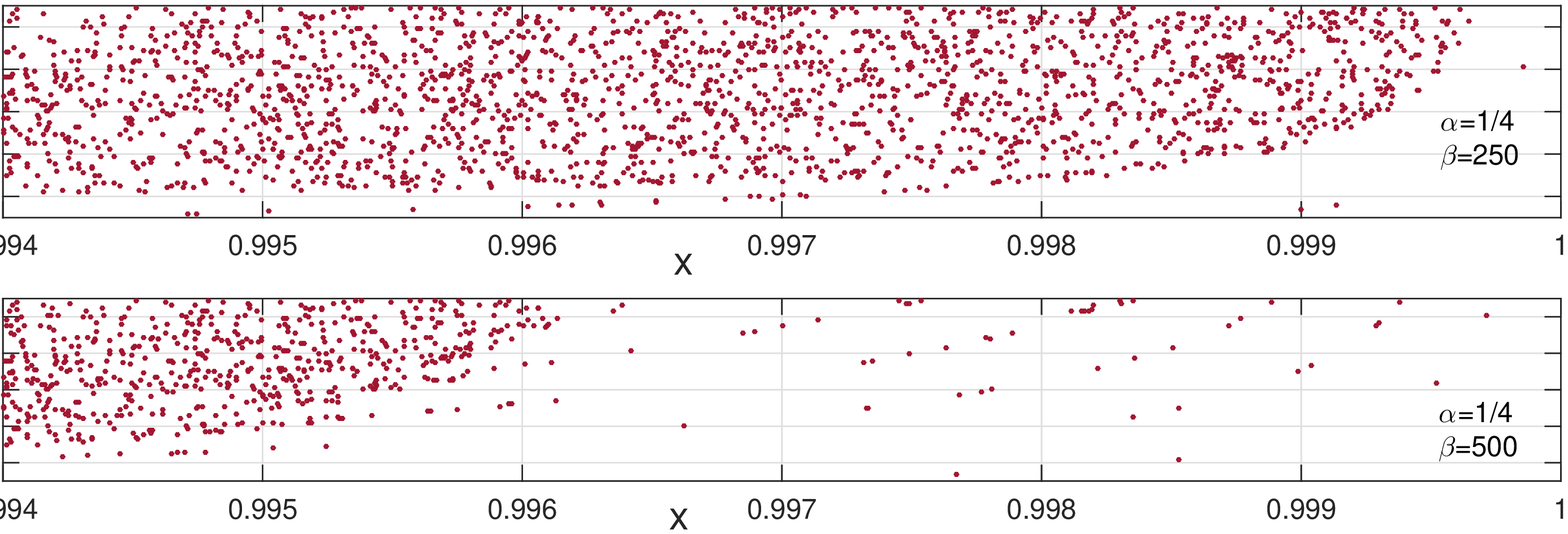}
\caption{
\label{fig:fig05} Performance of the asymptotic expansion (\ref{eq:Jacbeta1})
 for a fixed value of $\alpha$ and two different (large) values of the parameter $\beta$.  The variable $x$ 
 in the abscissas axis
is the argument $1-2z/b$
in the expansion. The points where the
error in the computation is greater than $10^{-8}$ are plotted. }
%\end{center}
\end{figure}

When the parameter values $\alpha$ and/or $\beta$ are large, it is interesting to consider the expansions given
in Sections~\ref{sec:jaclargeb} and \ref{sec:jaclargeab}.  An example of the performance of the expansion for
$P_n^{(\alpha,\beta)}(x)$ given in
(\ref{eq:Jacbeta1}) is shown in Figure~\ref{fig:fig05} for a fixed value of $\alpha$ and two different values of the parameter $\beta$
($\beta=250$, $\beta=500$.) The first five coefficients of the expansions (\ref{eq:Jacbeta2})
have been considered in the computations and the algorithm given in \cite{Gil:2017:ELP} is
used to compute the Laguerre polynomials.  The variable $x$ 
 in the abscissas axis
is the argument $1-2z/b$
in the expansion.
Random points have been generated   in the parameter region $(n,\,x) \in (10,\,110) \times (0.994,\,1)$.
The points where the
error in the computation is greater than $10^{-8}$ are plotted in the figure. 
 The results show that the range of validity of the expansion is quite narrow; therefore, it is of limited interest as a tool
to compute the polynomials.
The expansion is, however, very useful for obtaining expansions for the zeros of
Jacobi polynomials, as discussed in \cite{Gil:2018:ELB}.

\begin{figure}
%\begin{center}
\epsfxsize=14cm \epsfbox{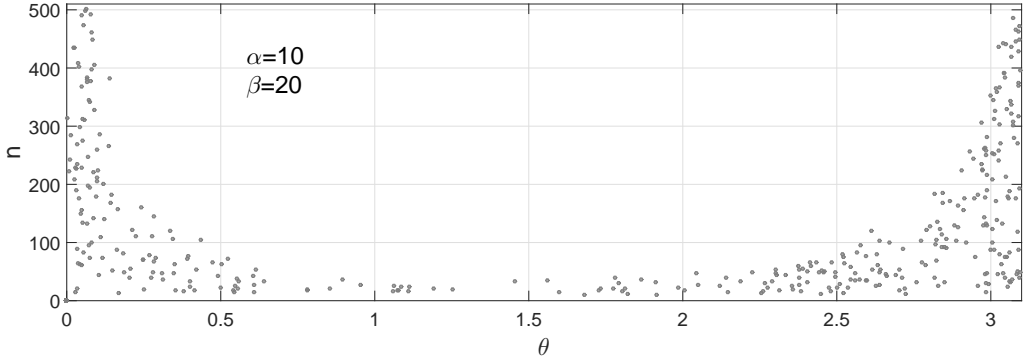}
\caption{
\label{fig:fig06} Example of the performance of the asymptotic expansion (\ref{eq:Jacalphabeta1}).
The points where the
error in the computation is greater than $10^{-8}$ are plotted. 
  }
%\end{center}
\end{figure}

More interesting for computational purposes is the expansion in terms of elementary functions  given in (\ref{eq:Jacalphabeta1}).
In Figure~\ref{fig:fig06} we show an example of the performance of this expansion. The first four coefficients of the expansions (\ref{eq:Jacalphabeta4})
have been considered in the computations. The points where the
error in the computation is greater than $10^{-8}$ are plotted in the figure. We have chosen $\theta$ ($x=\cos(\theta)$) as
the variable in the abscissas axis.

 As can be seen, the expansion performs
well even for not too large values of $\alpha$ and $\beta$ ($\alpha=10$, $\beta=20$ in the example shown.) Also,
as expected, the values of $x=\cos (\theta)$ should not be close to $\pm 1$.

\newpage
\section*{Appendix}

The following methods of approximation were used in our algorithm
for computing the Bessel functions $J_{\nu}(x)$:

\begin{description}

\item{\bf Power series.}

When $z$ is small, we use the power series given in Eq.(10.2.2) of \cite[\S10.19(ii)]{Olver:2010:Bessel}:

\[\mathop{J_{\nu}\/}\nolimits\!\left(z\right)=(\tfrac{1}{2}z)^{\nu}\sum_{k=0}^{%
\infty}(-1)^{k}\frac{(\tfrac{1}{4}z^{2})^{k}}{k!\mathop{\Gamma\/}\nolimits\!%
\left(\nu+k+1\right)}.\]

\item{\bf Debye's asymptotic expansions.}

The expressions are given in  Eq.(10.19.3) and Eq.(10.19.6) of
\cite[\S10.19(ii)]{Olver:2010:Bessel}: 

When $\nu <z$, we use

\[\mathop{J_{\nu}\/}\nolimits\!\left(\nu\mathop{\mathrm{sech}\/}\nolimits\alpha%
\right)\sim\frac{e^{\nu(\mathop{\tanh\/}\nolimits\alpha-\alpha)}}{(2\pi\nu%
\mathop{\tanh\/}\nolimits\alpha)^{\frac{1}{2}}}\sum_{k=0}^{\infty}\frac{U_{k}(%
\mathop{\coth\/}\nolimits\alpha)}{\nu^{k}},\]
and for $\nu >z$ 

\[\mathop{J_{\nu}\/}\nolimits\!\left(\nu\mathop{\sec\/}\nolimits\beta\right)\sim%
\left(\frac{2}{\pi\nu\mathop{\tan\/}\nolimits\beta}\right)^{\frac{1}{2}}\*%
\left(\mathop{\cos\/}\nolimits\xi\sum_{k=0}^{\infty}\frac{U_{2k}(i\mathop{\cot%
\/}\nolimits\beta)}{\nu^{2k}}-i\mathop{\sin\/}\nolimits\xi\sum_{k=0}^{\infty}%
\frac{U_{2k+1}(i\mathop{\cot\/}\nolimits\beta)}{\nu^{2k+1}}\right).\]

The coefficients $U_k(p)$ are polynomials in $p$ of degree $3k$ given by $U_0(p)=1$ and

\[U_{k+1}(p)=\tfrac{1}{2}p^{2}(1-p^{2})U_{k}^{\prime}(p)+\frac{1}{8}\int_{0}^{p}%
(1-5t^{2})U_{k}(t)dt.\]

\item{\bf Hankel's expansion.}

 For large values of the argument $z$, we use the Hankel's expansion given
in  \cite[\S10.17(i)]{Olver:2010:Bessel}:

\[\mathop{J_{\nu}\/}\nolimits\!\left(z\right)\sim\left(\frac{2}{\pi z}\right)^{%
\frac{1}{2}}\*\left(\mathop{\cos\/}\nolimits\omega\sum_{k=0}^{\infty}(-1)^{k}%
\frac{a_{2k}(\nu)}{z^{2k}}-\mathop{\sin\/}\nolimits\omega\sum_{k=0}^{\infty}(-%
1)^{k}\frac{a_{2k+1}(\nu)}{z^{2k+1}}\right),\]
where

\[\omega=z-\tfrac{1}{2}\nu\pi-\tfrac{1}{4}\pi.\]

The coefficients $a_{k}(\nu)$ are given by $a_0(\nu)=1$ and

\[a_{k}(\nu)=\frac{(4\nu^{2}-1^{2})(4\nu^{2}-3^{2})\cdots(4\nu^{2}-(2k-1)^{2})}{%
k!8^{k}},\quad k=1,2,3,\ldots\,.
\]

\item{\bf Airy-type expansions.}

We use the representation given in \cite[Chapter~8]{Gil:2007:NSF} 

$$
\dsp{J_\nu(\nu x)
=\quad \frac{\phi(\zeta)}{\nu^{1/3}}\left[\Ai(\nu^{2/3}\zeta)\,{A}_{\nu}(\zeta)
+\nu^{-4/3}\Ai'(\nu^{2/3}\zeta)\,{B}_{\nu}(\zeta)\right],}\\ 
\renewcommand{\arraystretch}{1.0}
$$

where 

$$
\phi(\zeta)=\left(\frac{4\zeta}{1-x^2}\right)^{\frac14},\quad \phi(0)=2^{\frac13}.
$$

The variable $\zeta$ is written in terms of the variable $x$ as

\renewcommand{\arraystretch}{1.5}
$$
\begin{array}{llll}
\dsp{\tfrac23\zeta^{3/2}}&=&\dsp{\ln\frac{1+\sqrt{{1-x^2}}}x-\sqrt{{1-x^2}},} &  \quad 0 < x\le1,\\
\dsp{\tfrac23(-\zeta)^{3/2}}&=&\dsp{\sqrt{{x^2-1}}-\arccos\frac{1}{x},} & \quad x\ge1.
\end{array}
\renewcommand{\arraystretch}{1.0}
$$

\item{\bf Three-term recurrence relation using Miller's algorithm.}

The standard three-term recurrence relation for the cylinder functions  

\[\mathop{{J}_{\nu-1}\/}\nolimits\!\left(z\right)+\mathop{{J}_{%
\nu+1}\/}\nolimits\!\left(z\right)=(2\nu/z)\mathop{{J}_{\nu}\/}%
\nolimits\!\left(z\right),\]

is computed backwards (starting from large values of $\nu$) using Miller's algorithm; see  \cite[\S4.6]{Gil:2007:NSF}.

\end{description}

\section{Acknowledgements}

We thank the referee for valuable comments. We acknowledge financial support from Ministerio de Ciencia e Innovaci\'on, Spain, 
projects MTM2015-67142-P (MINECO/FEDER, UE) and PGC2018-098279-B-I00 (MCIU/AEI/FEDER, UE). 
NMT thanks CWI, Amsterdam, for scientific support.

\bibliographystyle{plain}
\bibliography{gauss}

\end{document}